\documentclass[12pt]{amsart}
\usepackage{amsmath, amsthm, amssymb, amsopn, amsfonts, enumerate}
\usepackage{stmaryrd}
\newtheorem{thm}{Theorem}
\newtheorem{prop}[thm]{Proposition}

\newtheorem{conj}{Conjecture}

\newcommand{\Q}{\mathbb{Q}}

\newcommand{\Z}{\mathbb{Z}}
\newcommand{\M}{\mathcal{M}}
\renewcommand{\S}{\mathcal{S}}
\newcommand{\R}{\mathcal{R}}
\newcommand{\Mbar}{\overline{\mathcal{M}}}
\DeclareMathOperator{\FZ}{FZ}
\DeclareMathOperator{\Aut}{Aut}

\title{\textbf{Conjectural relations in the tautological ring of }$\Mbar_{g,n}$}
\author{A. Pixton}
\date{July 2012}
\begin{document}

\maketitle


\setcounter{section}{-1}
\section{Introduction}

The Faber-Zagier (FZ) relations in the tautological ring of the moduli space $\M_g$ of smooth curves of genus $g$ were recently proven to be true relations by R. Pandharipande and the author (see \cite{Berlin} for a sketch of the proof). These relations were constructed using the geometry of the moduli of stable quotients \cite{Marian-Oprea-Pandharipande}, which actually produces relations in $\Mbar_g$. Although we only analyzed these relations after restriction to $\M_g$, this means that our approach to the FZ relations also proves that these relations must extend tautologically to $\Mbar_g$.

These notes give a conjectural description of this extension and also explain how to add marked points to the relations. The result is a very large class of conjectural relations in the tautological ring of $\Mbar_{g,n}$. These notes are loosely based on informal talks given by the author at the workshop at KTH Stockholm on ``The moduli space of curves and its intersection theory'' in April 2012.

\subsection{Acknowledgements}
I would like to thank R. Pandharipande for suggesting the use of a formal strata algebra and for many helpful discussions about the tautological ring.

My research was supported by the Department of Defense through a NDSEG fellowship.

\section{The strata algebra}\label{taut classes}

We begin by reviewing a set of additive generators for the tautological ring $R^*(\Mbar_{g,n})$ that were described by Graber and Pandharipande (\cite{Graber-Pandharipande}, Proposition 11).

The pure boundary strata in $\Mbar_{g,n}$ are parametrized by their dual graphs: replace each irreducible component with a vertex labeled with the genus of the component, replace nodes with edges, and replace markings with half-edges labeled with the marking number. The only combinatorial constraints on such a dual graph $\Gamma$ with $n$ half-edges are that it should be connected, that any vertex of genus zero should have degree at least three, and that the sum of the genera of the vertices plus the cycle number of the graph should be equal to $g$. We will let $V(\Gamma)$ denote the set of vertices and $E(\Gamma)$ the set of full edges of $\Gamma$.

A valid dual graph $\Gamma$ gives a gluing map $\xi_\Gamma:\Mbar_\Gamma := \Mbar_{g_1,n_1}\times \Mbar_{g_2,n_2}\times \cdots \times \Mbar_{g_k,n_k} \to \Mbar_{g,n}$. The set of generators given in \cite{Graber-Pandharipande} consists of the pushforwards ${\xi_\Gamma}_*\theta$ of classes $\theta$ that are monomials in the Arbarello-Cornalba $\kappa$ (see \cite{Arbarello-Cornalba}) and cotangent line $\psi$ classes on the components of $\Mbar_\Gamma$. The automorphism group of $\Gamma$ acts on these monomials, and we want to take just one representative from each orbit under this action.

Let $\S_{g,n}$ be the finite-dimensional $\Q$-vector space with basis labeled by the generators described above. The proof in \cite{Graber-Pandharipande} that the linear span of these generators is closed under multiplication in the tautological ring provides rules for multiplying two of the generators together and re-expressing as a linear combination of the generators, and we can use this to define a $\Q$-bilinear multiplication operation on $\S_{g,n}$. It is straightforward to verify that this operation is associative (by defining $(A,B,C)$-graphs analogously to the $(A,B)$-graphs in \cite{Graber-Pandharipande}). Thus $\S_{g,n}$ with this multiplication has the structure of a commutative $\Q$-algebra. We call $\S_{g,n}$ the {\it strata algebra} and note that it comes equipped with a natural surjection to the tautological ring $R^*(\Mbar_{g,n})$.

In addition, the rules given in \cite{Graber-Pandharipande} and \cite{Arbarello-Cornalba} for taking the pushforward or pullback of one of the generators along a forgetful or gluing map allow us to define $\Q$-linear maps between the strata algebras lifting the pushforward and pullback maps between the tautological rings. It is again straightforward to check that these pushforward and pullback maps on strata algebras satisfy the basic ring-theoretic properties one would expect (i.e. pullbacks are ring homomorphisms and the projection formula holds).

These strata algebras should be viewed as generalizations to $\Mbar_{g,n}$ of the formal polynomial algebra $\Q[\kappa_1,\kappa_2,\ldots]$ that surjects onto the tautological ring of $\M_g$. From now on, by a {\it tautological relation} we will mean an element of the kernel of the natural surjection $\S_{g,n}\to R^*(\Mbar_{g,n})$.

\section{The relations}\label{relations}
\subsection{The Faber-Zagier (FZ) relations}
Before describing our generalization of the FZ relations, we recall how the FZ relations themselves are constructed. The following description can be easily shown to be equivalent to the usual ones, though it may appear slightly different.

We need to establish some notation. First, $A$ and $B$ are the fundamental power series
\[
A = \sum_{n\ge 0}\frac{(6n)!}{(3n)!(2n)!}T^n
\]
and
\[
B = \sum_{n\ge 0}\frac{6n+1}{6n-1}\cdot\frac{(6n)!}{(3n)!(2n)!}T^n.
\]
Then $(C_i)$ is a sequence of power series for $i\ge 0$, $i\not\equiv 2$ mod $3$, defined in terms of $A$ and $B$:
\[
C_{3i} = T^iA,
\]
\[
C_{3i+1} = T^iB.
\]

Next, if $F$ is a power series in $T$ then we let $[F]_{T^n}$ denote the coefficient of $T^n$ and define
\[
\{F\} = \sum_{n}[F]_{T^n}K_nT^n,
\]
where the $K_n$ are formal indeterminates.

Finally, $\kappa$ is a linear operator converting polynomials in $K_0,K_1,\ldots$ into polynomials in $\kappa_0,\kappa_1,\ldots$, defined by
\[
\kappa(K_{e_1}\cdots K_{e_l}) = \sum_{\tau \in S_l}\prod_{c\text{ cycle in }\tau}\kappa_{e_c}
\]
where $e_c$ is the sum of the $e_i$ appearing in the cycle $c$ in a permutation of $\{1,\ldots,l\}$.

The FZ relations in $R^r(\M_g)$ are parametrized by partitions $\sigma$ with no parts of size $2$ mod $3$, subject to the conditions $3r\ge g+1+|\sigma|$ and $3r\equiv g+1+|\sigma|$ mod $2$. Given such a partition $\sigma$ with parts $\sigma_1,\ldots,\sigma_l$, the corresponding FZ relation is
\[
\FZ(g,r; \sigma) = \bigg[\kappa\Big(\exp\big(\{1-A\}\big)\{C_{\sigma_1}\}\cdots\{C_{\sigma_l}\}\Big)\bigg]_{T^r}.
\]

\subsection{Conjectural relations}
We now begin to construct conjectural relations in $R^r(\Mbar_{g,n})$. When $n=0$, it will be clear that the restriction of these relations to the interior are simply the usual FZ relations described above. The relations are parametrized by partitions $\sigma$ with no parts of size $2$ mod $3$ together with nonnegative integers $a_1,\ldots,a_n$ not $2$ mod $3$, subject to the conditions $3r\ge g+1+|\sigma|+\sum_i a_i$ and $3r\equiv g+1+|\sigma|+\sum_i a_i$ mod $2$.

We denote the relation coming from this data by
\[
\R(g,n,r; \sigma, a_1, \ldots, a_n),
\]
and we write
\[
\R(g,n,r; \sigma, a_1,\ldots,a_n) = \sum_\Gamma\frac{1}{|\Aut(\Gamma)|}{\xi_\Gamma}_*\big(\R_\Gamma(g,n,r; \sigma, a_1,\ldots,a_n)\big),
\]
where the sum is over isomorphism classes of dual graphs $\Gamma$. Here $\R_\Gamma(g,n,r; \sigma, a_1,\ldots,a_n)$ is a polynomial of degree $r - |E(\Gamma)|$ in the $\kappa$ and $\psi$ classes on the components of $\Mbar_\Gamma$. We denote these classes by $\kappa_i^{(v)}$ (for $i\ge 0$ and $v\in V(\Gamma)$ a vertex of $\Gamma$) and $\psi_h$ (for $h$ a half-edge of $\Gamma$, by which we mean either half of a full edge or one of the half-edges $h_i$ corresponding to a marking $i$).

In order to describe this polynomial, we need to modify the definitions used in the FZ relations to include extra parity information. This will involve augmenting the polynomials and power series with extra commuting variables $\zeta$ satisfying $\zeta^2 = 1$.

First, we define series $\widehat{C}_i(T,\zeta)\in \left(\Q[\zeta]/(\zeta^2 - 1)\right)[[T]]$ by
\[
\widehat{C}_{3i}(T,\zeta) = T^iA(\zeta T)
\]
and
\[
\widehat{C}_{3i+1}(T,\zeta) = \zeta T^iB(\zeta T).
\]

Next, if $F$ is a power series in $T$ and $\zeta$ then we let $[F]_{T^n\zeta^a}$ denote the coefficient of $T^n\zeta^a$ and define
\[
\{F\} = \sum_{n\in\Z, a\in\Z/2}[F]_{T^n\zeta^a}K_{n,a}T^n,
\]
where the $K_{n,a}$ are formal indeterminates.

Finally, $\widehat{\kappa}$ is a linear operator converting polynomials in the $K_{n,a}$ into polynomials in the kappa variables $\kappa_i^{(v)}$ along with an additional variable $\zeta_v$ (satisfying $\zeta_v^2 = 1$) for each vertex $v$, defined by
\[
\widehat{\kappa}(K_{e_1,a_1}\cdots K_{e_l,a_l}) = \sum_{\tau \in S_l}\prod_{c\text{ cycle in }\tau}\left(\sum_{v\in V(\Gamma)}\kappa_{e_c}^{(v)}\zeta_v^{a_c}\right)
\]
where $e_c$ and $a_c$ are the sums of the $e_i$ and $a_i$ respectively appearing in the cycle $c$.

Then we can write
\begin{align*}
\R_\Gamma(g,n,r; \sigma, a_1,&\ldots,a_n) = \bigg[\frac{1}{2^{h_1(\Gamma)}}\widehat{\kappa}\Big(\exp\big(\{1 - \widehat{C}_0\}\big)\{\widehat{C}_{\sigma_1}\}\cdots\{\widehat{C}_{\sigma_l}\}\Big)\\ &\cdot\prod_{i = 1}^n\widehat{C}_{a_i}(\psi_{h_i}T,\zeta_{v_i})\prod_{e\in E(\Gamma)}\Delta_e    \bigg]_{\displaystyle T^{r - |E(\Gamma)|}\prod_{v\in V(\Gamma)}\zeta_v^{g_v+1}}\quad,
\end{align*}
where $h_1(\Gamma) = |E(\Gamma)| - |V(\Gamma)| + 1$ is the cycle number of $\Gamma$, marking $i$ corresponds to half-edge $h_i$ on vertex $v_i$, and $g_v$ is the genus of vertex $v$.

Also, for each edge $e \in E(\Gamma)$, let $e_1$ and $e_2$ be the two halves, attached to vertices $v_1$ and $v_2$ respectively. The edge contribution $\Delta_e$ appearing in the above formula is a power series in $T$ with coefficients that are polynomials in $\psi_1 := \psi_{e_1}, \psi_2 := \psi_{e_2}, \zeta_1 := \zeta_{v_1}$, and $\zeta_2 := \zeta_{v_2}$:
\[
\Delta_e = \frac{A(\zeta_1\psi_1T)\zeta_2B(\zeta_2\psi_2T) + \zeta_1B(\zeta_1\psi_1T)A(\zeta_2\psi_2T) + \zeta_1 + \zeta_2}{(\psi_1 + \psi_2)T}.
\]
The fact that $\psi_1 + \psi_2$ divides the numerator in this formula is a consequence of the identity
\[
A(T)B(-T) + A(-T)B(T) + 2 = 0.
\]

This completes the definition of $\R(g,n, r; \sigma, a_1, \ldots, a_n)$.

\begin{conj}\label{true}
$\R(g,n, r; \sigma, a_1, \ldots, a_n)$ maps to $0 \in R^*(\Mbar_{g,n})$ if $3r\ge g+1+|\sigma|+\sum_i a_i$. 
\end{conj}
In principle, it seems that it should be possible to adapt the stable quotients methods developed in \cite{Marian-Oprea-Pandharipande} and used in \cite{Berlin} to prove this, though there are a lot of details to work out.

We also conjecture that this construction produces all relations in the tautological ring. To make the statement of this precise, let $\R_{g,n}$ be the linear span of all elements of the strata algebra $\S_{g,n}$ produced as follows: choose a dual graph $\Gamma$ for a boundary stratum of $\Mbar_{g,n}$, pick one of the components $\S_{g',n'}$ in $\S_\Gamma = \S_{g_1,n_1}\times \cdots \times \S_{g_m,n_m}$, take the product of a relation $\R(g',n', r; \sigma, a_1, \ldots, a_n)$ on the chosen component together with arbitrary classes on the other components, and push forward along the gluing map $\S_\Gamma \to \S_{g,n}$.
\begin{conj}\label{all}
$\R_{g,n}$ is the kernel of the natural surjection $\S_{g,n}\to R^*(\Mbar_{g,n})$.
\end{conj}

Note that this conjecture would imply that the regular FZ relations give all relations in the tautological ring of $\M_g$, which is well known to contradict for $g\ge 24$ Faber's conjecture \cite{Faber} that this ring is Gorenstein. Another consequence is that all ``new'' relations in positive genus would be $S_n$-invariant - see Proposition~\ref{new}.

\section{Additional properties}\label{properties}
In this section we list a couple of properties satisfied by the relations $\R(g,n, r; \sigma, a_1, \ldots, a_n)$.

First, there are no tautological ways of enlarging the space of conjectured relations $\R_{g,n}$.
\begin{prop}
The vector subspace $\R_{g,n}$ of the strata algebra $\S_{g,n}$ is an ideal. Moreover, this collection of ideals of the stata algebra is closed under pushforward and pullback by the gluing and forgetful maps.
\end{prop}

We can also ask how many of these relations are actually necessary to generate them all using multiplication, pushforward, and pullback. To make this question precise, let $\R^{\text{old}}_{g,n}$ be the sub-ideal of $\R_{g,n}$ generated by $\S^{>0}_{g,n}\R_{g,n}$ together with the images of other $\R_{g',n'}$ under pushforwards via gluing maps or pullbacks via forgetful maps; these are the relations that come from some simpler moduli space.
\begin{prop}\label{new}
If $g>0$ then $\R_{g,n}/\R^{\text{old}}_{g,n}$ is generated by the relations $\R(g,n,r; \sigma, 1,\ldots, 1)$ with all parts of $\sigma$ congruent to $1$ mod $3$.
\end{prop}

In particular, all the new relations in our set of conjectural relations are $S_n$-invariant.

\section{Computations}\label{computations}
We have checked both Conjectures~\ref{true} and \ref{all} in many small cases (where the tautological ring has already been computed) with the aid of code written for Sage to compute the finitely many relations for any fixed $g,n$, and $r$. This includes checking for the presence of Getzler's relation \cite{Getzler} in $R^2(\Mbar_{1,4})$ and the Belorousski-Pandharipande relation \cite{Belorousski-Pandharipande} in $R^2(\Mbar_{2,3})$. In addition, Yang computed the ranks of the Gorenstein quotient in many cases in \cite{Yang} and in most of them we have been able to compute that the relations $\R_{g,n}$ produce the same ranks. For example, the relations $\R_{g,n}$ give rank $333$ for $R^3(\Mbar_{2,4})$, rank $142$ for $R^4(\Mbar_{3,2})$, rank $50$ for $R^4(\Mbar_4)$ and $R^5(\Mbar_4)$, and so on.

\bibliographystyle{amsplain}
\bibliography{tautrel}

\vspace{+8 pt}
\noindent
Department of Mathematics\\
Princeton University\\
apixton@math.princeton.edu

\end{document}